\documentclass[nato,numreferences]{crckbked}
 \usepackage{amsmath}
\usepackage{amsfonts}
\usepackage{amsthm}
\usepackage{amstext}

\newtheorem{definition}{Definition}
\newtheorem{theorem}{Theorem}
\newtheorem{proposition}{Proposition}
\theoremstyle{remark}

\newtheorem{example}{\sc Example}

\newcommand{\hb}{\ensuremath{\mathcal{H}}}
\newcommand{\tenz}{\ensuremath{T_{ij}^{kl}}}

\begin{document}

\begin{article}
\begin{opening}

\title{A family of $*$-algebras allowing Wick ordering: Fock representations
and universal enveloping $C^*$-algebras}

\author{Palle E. T. \surname{Jorgensen}
\thanks{\email{jorgen@math.uiowa.edu}}}
\institute{Department of Mathematics\\The University of Iowa\\
Iowa City, Iowa 52242-1419 U.S.A.}
\author{Daniil P. \surname{Proskurin}
\thanks{\email{prosk@imath.kiev.ua}}}
\institute{Kyiv Taras Shevchenko University\\
Cybernetics Department\\Volodymyrska, 64, Kyiv, 01033, Ukraine}
\author{Yurii \surname{Samoilenko}
\thanks{\email{yurii\_sam@imath.kiev.ua}}}
\institute{Institute of Mathematics, National Academy of Sciences\\
Tereschenkivska, 
3, Kyiv, 01601, Ukraine}

\runningauthor{Palle E.T. Jorgensen, Daniil P. Proskurin, Yurii Samoilenko}
\runningtitle{A family of $*$-algebras allowing Wick ordering}

\end{opening}

\begin{quotation}
\textsc{Abstract.}
We consider an abstract Wick
ordering as a family of relations on elements $a_{i}$ and
define $*$-algebras by these relations. The relations
are given by a fixed operator $T\colon
\mathfrak{h}\otimes\mathfrak{h}\rightarrow
\mathfrak{h}\otimes\mathfrak{h}$, where $\mathfrak{h}$ is one-particle
space, and they naturally define both a $*$-algebra and an
inner-product space $\mathcal{H}_{T}$,
$\langle\,\cdot\,,\,\cdot\,\rangle_{T}$.
If $a_{i}^{*}$ denotes
the adjoint, i.e.,
$\langle a_{i}\varphi,\psi\rangle_{T}=
\langle\varphi,a_{i}^{*}\psi\rangle_{T}$,
then we identify when
$\langle\,\cdot\,,\,\cdot\,\rangle_{T}$
is positive semidefinite (the positivity question!).
In the case of deformations of the CCR-relations
(the $q_{ij}$-CCR and the twisted CCR's), we
work out the universal
$C^{*}$-algebras $\mathfrak{A}$, and we prove that, in
these cases, the Fock representations of the
$\mathfrak{A}$'s are faithful.
\end{quotation}

\section{Introduction}\label{jor-Int}
In recent papers
\cite{jor-bied,jor-burkl,jor-fiv,jor-jps,jor-mcf,jor-pw},
the applications of Lie superalgebras, quantum
groups, $q$-algebras in mathematical physics have stimulated
interest in the $*$-algebras defined by generators and relations and
their representations by Hilbert space operators. For example, the
representations of various deformations of canonical commutation
relations (CCR), in particular Fock representaion, were used to
construct non-classical models of theoretical physics and probability,
such as the free quon gas (see \cite{jor-rfw1}), $q$-Gaussian processes
(see \cite{jor-bsp}) etc.

The constructions are
interesting from both physical and mathematical points of view.
They give a
canonical realisation of a given deformed relation like the
Fock representation, or a realisation by differential operators. When the
relations can be realised by bounded operators, it
is useful to study the universal enveloping $C^*$-algebras for them
and the stability of isomorphism classes of these $C^*$-algebras on
parameters (see for example \cite{jor-dn,jor-jsw2}).
The stability question \cite{jor-jsw2} refers to how the $C^{*}$-isomorphism
classes depend on variations in the deformation
variables; in some cases there are open
regions in parameter space where the $C^{*}$-isomorphism
class is constant.

In the present paper we give a review of some results concerning a wide
class of deformed relations of the following form
\begin{equation}
a_i^* a_j = \delta_{ij}1+\sum_{k,l=1}^d T_{ij}^{kl} a_la_k^*,\
i,j=1,\ldots,d,\label{ijor-eqj1}
\end{equation}
where  $T_{ij}^{kl}\in \mathbb{C}$,  such that $T_{ij}^{kl} =
\bar{T}_{ji}^{lk}$. These relations generate a $*$-algebra allowing Wick
ordering or {\bf Wick algebra} (see \cite{jor-jps,jor-jsw,jor-bs,jor-mral}).
The $*$-algebra $\mathfrak{A}_{T}$ has a naturally defined
Fock vacuum ``state'' or functional and there is a
corresponding inner-product space $\mathcal{H}_{T}$,
$\langle\,\cdot\,,\,\cdot\,\rangle_{T}$,
such that, in the associated GNS-representation,
the identity
$\langle a_{i}\varphi,\psi\rangle_{T}=
\langle\varphi,a_{i}^{*}\psi\rangle_{T}$
holds.
But the vacuum
functional is generally not positive, and
the operators in the representation not bounded, and
therefore the Hermitian inner product
$\langle\,\cdot\,,\,\cdot\,\rangle_{T}$
is then generally not positive semidefinite. The
positivity question, and the faithfulness of the
Fock representation, are the foci of this paper.

Note that (\ref{jor-eqj1}) generalises some well-known types
of deformed commutation relations, quantum groups, etc.
(see \cite{jor-bied,jor-fiv,jor-mcf,jor-pw,jor-bsp,jor-bs,jor-green,jor-marc}).
The basic examples for
us will be the $q_{ij}$-CCR introduced and studied by
M. Bo\.zejko and R. Speicher (see \cite{jor-bsp,jor-bs}), and the twisted
canonical commutation relations (TCCR) constructed by W. Pusz and
S.L. Woronowicz (see \cite{jor-pw}).
They were further studied in \cite{jor-jswc}
where the traditional Cuntz algebra of
\cite{jor-cun} was considered as a base-point,
corresponding to $q_{ij}=0$,
and the variation of the $C^*$-isomorphism
class was considered as a function of $q_{ij}$.

\begin {example}
\label{jor-Exa1}
$q_{ij}$-CCR, $2d$ generators:
\begin{align*}
\mathbb{C}\langle &a_i,
\ a_i^*\mid a_i^*a_j=\delta_{ij}1+q_{ij}a_ja_i^*,\ i,j=1,\ldots,d,\\
& q_{ji}=\overline{q}_{ij}\in\mathbb{C},\ \left\vert q_{ij}\right\vert \le
 1\rangle
\end{align*}
\end{example}
\begin{example}
The Wick algebra for TCCR:
\begin{align*}\label{jor-mucc}
a_i^*a_i &=1+\mu^2 a_ia_i^* -(1-\mu^2)\sum_{k<i}a_k a_k^*,\
i=1,\ldots,d  \\
a_i^*a_j &=\mu a_ja_i^*,\ i\ne j,\quad
 0<\mu<1
\end{align*}
\end{example}

We present some sufficient conditions on the coefficients
$\{T_{ij}^{kl}\}$ for the existence of the Fock representation, and we
describe the structure of the Fock space. We also give conditions for
the faithfulness of Fock representation and describe its kernel in the
degenerated case (see Sec.\ \ref{jor-fock}).

Further we consider the universal $C^*$-algebras for
the examples above.
Specifically we show
that the universal $C^*$-algebras for
$q_{ij}$-CCR (TCCR) can
be generated by isometries (partial isometries) satisfying
a certain algebraic relation.
The description of the
$C^*$-isomorphism classes for different values
of parameters is presented.

We also show that the Fock representations of
$q_{ij}$-CCR for some values of parameters, and TCCR for any value of
parameter, are faithful on the $C^*$-level, i.e., the Fock
representations of the corresponding $C^*$-algebras are faithful
(see Sec.\ \ref{jor-stab}).

The complete proofs of all results presented here can be found in
\cite{jor-jps,jor-jsw2,jor-jsw,jor-dpl,jor-prsm}. For detailed information
about $*$-representations of finitely
generated $*$-algebras see \cite{jor-osam1}.

\section{Basic definitions}\label{jor-Bas}
Firstly let us construct a canonical realisation of Wick algebra,
i.e., the $*$-algebra on the relations (\ref{jor-eqj1}), with
coefficents $\{T_{ij}^{kl}\}$: we denote it by $W(T)$. To do it
consider a finite-dimensional Hilbert space
$\mathcal{H}=\langle e_{1}, \ldots ,e_{d} \rangle  $. Construct the full tensor
algebra over $\mathcal{H},\ \mathcal{H}^*$, denoted by
$\mathcal{T}(\mathcal{H},\mathcal{H}^{*})$. Then
    \begin{equation}
W(T) \cong
\mathcal{T}(\mathcal{H},\mathcal{H}^{*}) /
\langle e_{i}^{*}\otimes e_{j} -
\delta_{ij}1 -\sum \tenz e_{i}\otimes e_{j}^{*} \rangle ,\label{ijor-eqj2}
    \end{equation}
dividing out by the two-sided ideal on the relations (\ref{jor-eqj1}).
Note that in this realisation the subalgebra of $W(T)$ generated
by $\{a_i\}$ is identified with the $\mathcal{T}(\mathcal{H})$.

The following operators were presented in \cite{jor-jsw} as a useful
tool for computation with Wick algebras and their Fock representations.
\begin{align}
T\colon & \mathcal{H} \otimes \mathcal{H}\mapsto \mathcal{H}
 \otimes\mathcal{H}, \quad  T e_{k}\otimes e_{l} =
\sum_{i,j} T_{ik}^{lj}e_{i}\otimes e_{j},\ T=T^*\nonumber
\\
T_{i}\colon & \mathcal{H}^{\otimes n}\mapsto\mathcal{H}^{\otimes n},
\quad
T_{i}=\underbrace{1\otimes \cdots\otimes 1}_{i-1}\otimes T
\otimes\underbrace{1\otimes\cdots\otimes 1}_{n-i-1},\nonumber
\\
R_{n}\colon & \mathcal{H}^{\otimes n}\mapsto \mathcal{H}^{\otimes n},
\quad
R_{n}=1+T_1+T_1 T_2+\cdots + T_1 T_2\cdots T_{n-1},\nonumber
\\
P_n\colon & \mathcal{H}^{\otimes n}\mapsto \mathcal{H}^{\otimes n},
\quad
P_2=R_2,\ P_{n+1}=(1\otimes P_n)R_{n+1}.\label{ijor-eqj3}
\end{align}
The sequences of operators $P_{0}=1_{\operatorname*{vac}}$,
$R_{1}=1+T$, $P_{1}=(1\otimes 1)(1+T)\cong1+T$,
$R_{2}$, $\dots$, $P_{n}$
are defined recursively. It is the sequence $P_{n}$
which enters into the positivity question. The other
one is only intermediate. The Hermitian inner
product
$\langle\,\cdot\,,\,\cdot\,\rangle_{T}$
on $\mathcal{T}_{n}(\mathcal{H})$ is then
\begin{equation}
\langle\phi,\psi\rangle_{\mathcal{T}_{n}(\mathcal{H})}
:=\langle\varphi,P_{n}\psi\rangle_{\operatorname*{tensor}}\label{ijor-eqj4}
\end{equation}
where
$\langle\,\cdot\,,\,\cdot\,\rangle_{\operatorname*{tensor}}$
is just the usual inner product on
$\mathcal{T}_{n}(\mathcal{H})$ induced by
$\langle\,\cdot\,,\,\cdot\,\rangle$
on $\mathcal{H}$.
Hence, we need conditions on $T\colon
\mathcal{H}\otimes\mathcal{H}\rightarrow\mathcal{H}\otimes\mathcal{H}$
which make the operators $P_{n}$ positive for all $n$.
For example, in terms of these operators we can describe the procedure
of Wick ordering, i.e., the commutation formula for fixed generator
$a_i^*$ and any homogeneous polynomial in $a_k,\ k=1,\ldots,d$ (see
\cite{jor-proid}).
\begin{proposition}
\label{jor-Pro1}
Let $X\in\mathcal{H}^{\otimes n}$. Then
\begin{equation} \label{ijor-commutation_rule}
e_i^*\otimes X = \mu (e_i^*) R_n X + \mu (e_i^*)\sum_{k=1}^{d}
T_1 T_2 \cdots T_n (X\otimes e_k) e_k^*,
\end{equation}
where $\mu (e_i^*)\colon\mathcal{T}(\mathcal{H})
\mapsto\mathcal{T}(\mathcal{H})$ is defined as follows
\[
\mu (e_i^*) 1=0,\quad \mu (e_i^*) e_{i_1}\otimes \cdots\otimes e_{i_n}=
\delta_{ii_1} e_{i_2}\otimes\cdots\otimes e_{i_n}.
\]
\end{proposition}
For our examples the operator $T$ have the following form:
\begin{example}
\[
T e_i\otimes e_j = q_{ij} e_j\otimes e_i, i,j=1,\ldots,d.
\]
\end{example}
\begin{example}
\begin{align*}
T e_i\otimes e_i &  = \mu^2 e_i\otimes e_i \\
T e_i\otimes e_j & = \mu e_j\otimes e_i\ ,\ i<j \\
T e_i\otimes e_j & = -(1-\mu^2) e_i\otimes e_j + \mu e_j\otimes e_i\
,\ i>j .
\end{align*}
\end{example}
Note that for both examples, the operator $T$ satisfies a {\bf braid
condition}, i.e., on the $\mathcal{H}^{\otimes 3}$ we have
\begin{equation}
T_1T_2T_1=T_2T_1T_2.\label{ijor-eqj5}
\end{equation}
The operators presented above appear naturally in construction of Fock
representation of $W(T)$. This notion is induced in the obvious way
from the classical one for CCR, however, in general, the Fock space is
not always symmetric (see \cite{jor-jsw}).
\begin{definition}
The representation $\lambda_0$ acting on the space $\mathcal{T}(\hb)$
by formulas
\begin{align*}
\lambda_0(a_i)e_{i_1}\otimes\cdots\otimes
e_{i_n}& =e_i\otimes e_{i_1}\otimes\cdots\otimes e_{i_n},\quad
n\in\mathbb{N}\cup\{0\}\\
\lambda_0(a_i^*)1_{\operatorname*{vac}}&= 0
\end{align*}
where the action of $\lambda_0(a_i^*)$ on the monomials of degree
$n\ge 1$ is determined inductively using the basic relations,
is called the Fock representation.
\end{definition}
It is easy to see that $\lambda_0(a_i)$ are the classical creation
operators and $\lambda_0(a_i^*)$ are twisted annihilation ones.
Evidently in this way we have constructed a representation of $W(T)$,
but not yet a $*$-representation. To do it one has to supply the
$\mathcal{T}(\mathcal{H})$ by the appropriate inner product (see
\cite{jor-jsw}).
This is where formula (\ref{jor-eqj4}) comes in.
\begin{definition}
The
Fock inner product \textup{(}see \textup{\cite{jor-jsw})}
is the unique semilinear Hermitian form
$\langle \ ,\ \rangle  _{T}$ on
$\mathcal{T}(\mathcal{H})$ such that
\[
\langle  \lambda_0 (a_i)X, Y\rangle _{T} =
\langle  X, \lambda_0 (a_i^*)Y\rangle _{T},\quad
X,Y\in\mathcal{T}(\mathcal{H}).
\]
\end{definition}
Similarly to the definition of Fock representation, the Fock inner
product on $\mathcal{T}(\mathcal{H})$ can be computed inductively.
It is easy to see that for $X\in\mathcal{H}^{\otimes m}$,
$Y\in\mathcal{H}^{\otimes n}$, $n\ne m$, we have $\langle X,Y\rangle _{T}=0$.
On the components of powers $0,1$, the Fock inner product concides
with the standard one. For any $X,Y\in\mathcal{H}^{\otimes n}$, $n\ge 2$,
we have
\[
\langle X,Y\rangle _{T} =\langle X, P_n Y\rangle ,
\]
which agrees with (\ref{jor-eqj4}) above.
The operator $P_{n}=P_{n}(T)$ are given in (\ref{jor-eqj3}).

Evidently, if we want to extend the Fock representation of $W(T)$ to
the $*$-representation by Hilbert-space operators, we should require
that all the operators $P_{n}$, $n=2,\dots$, be
positive semidefinite, and that
the subspace
\[
\mathcal{I}=\bigoplus_{n\ge 2}\ker P_n
\]
determines the kernel of the Fock inner product. Consequently the
Hilbert-space
structure of the Fock space emerges.

\section{The structure of the Fock representation}\label{jor-fock}

In this section we present some sufficient conditions posed on the
operator $T$ for the positive-definite
property
of the Fock inner product, and we
show that the kernel of the Fock representation is generated as a
$*$\textit{-ideal} by the kernel of
the Fock inner product. In particular, when the
Fock inner product is strictly positive definite
(i.e., when it has zero kernel),
the Fock representation $\pi_{F}$ is
faithful, i.e., $\ker(\pi_{F})=0$.

There are several sufficient conditions on the operator $T$ for the
Fock inner product to be positive.
It was shown in \cite{jor-jsw2} that for sufficiently small coefficients we have
strict positivity of the Fock inner product. This result is a corollary
of the stability of the universal enveloping $C^*$-algebra for the
Wick algebra around the zero base point (see Sec.\ \ref{jor-stab}).
\begin{theorem}\label{jor-cunst}
If the operator $T$ satisfies the norm bound $\Vert T\Vert<\sqrt{2}-1$,
then $P_n>0$, $n\ge 2$,
where $>$ refers to strict positivity.
\end{theorem}
Another kind of sufficient condition is positivity of operator $T$
 (see \cite{jor-jsw}).
\begin{theorem}
If $T\ge 0$ then $P_n>0$, $n\ge 2$.
\end{theorem}
In the present paper we will suppose that the
operator $T$ satisfies the
{\bf braid} condition
(\ref{jor-eqj5}).
It was shown by M. Bo{\.z}ejko
and R. Speicher (see \cite{jor-bs}) that, in this case, the
operators $P_n$, $n\ge 2$, have a natural description in terms
of quasimultiplicative operator-valued mappings on the
Coxeter group $S_n$. The following is a corollary of a much
more general result proved in \cite{jor-bs} for mappings on the
general Coxeter group.
\begin{theorem}\label{jor-fobs}
Let $T$ satisfy the braid condition
\textup{(\ref{jor-eqj5})}
and suppose $-1\le T\le 1$. Then
$P_n\ge 0$. Moreover,
if $\Vert T\Vert\le 1$, then $P_n>0$,
and the operators of the Fock representation are bounded,
i.e., the Fock representation is by bounded operators.
\textup{(}Recall, the Fock representation of the undeformed
CCR-algebra is unbounded.\textup{)}
\end{theorem}
We present a more precise version of this theorem. Namely, we give
the description of kernel of $P_n$ in the degenerate case. As an
immediate corollary of this result we have the strict
positivity of $P_n$, $n\ge 2$, for braided $T$ satisfying the inequality
$-1<T\le 1$ (see \cite{jor-jps}).
\begin{theorem}\label{jor-bas}
Let $W(T)$ be a Wick algebra with braided operator
$T$ satisfying the norm bound $\Vert T\Vert\le 1$.
Then for any $n\ge 1$,
\[
\ker P_{n+1}=\sum_{k+l=n-1}
\hb^{\otimes^k}\otimes \ker(1+T)\otimes\hb^{\otimes^l}=
\sum_{k=1}^{n} \ker (1+T_k).
\]
\end{theorem}
Let us illustrate this result on the examples.
\begin{example}
For $q_{ij}$-CCR we have the alternatives:
\begin{itemize}
\item $\vert q_{ij}\vert<1$ for any $i,j=1,\ldots,d$.\\
In this case $-1<T<1$ and the Fock inner product is strictly positive.
\item $\vert q_{ij}\vert =1$, $i\ne j$.\\
For these values of parameters we have $-1\le T\le 1$ and
\[
\ker (1+T)=\langle  a_ja_i-q_{ij}a_ia_j,\ i<j\rangle .
\]
\end{itemize}
\end{example}
\begin{example}
For the TCCR Wick algebra, we have $-1\le T\le 1$, and
\[
\ker(1+T)=\langle  a_ja_i-\mu a_ia_j,\ i<j\rangle .
\]
\end{example}
The following proposition shows that, for algebras with braided operator
$T$, the kernel of the Fock representation is generated as a $*$-ideal
by the kernel of the Fock inner product, i.e.,
\[
\mathcal{I}=\bigoplus_{n\ge 2} \ker P_n.
\]
\begin{proposition}
Let $W(T)$ be a Wick algebra with braided operator $T$ and let the Fock
representation $\lambda_0$ be positive
\textup{(}i.e., the Fock inner product is
positive definite\/\textup{)}. Then
\[
\ker \lambda_0 = \mathcal{I}\otimes\mathcal{T}(\hb^*) +
\mathcal{T}(\hb)\otimes\mathcal{I}^*.
\]
\end{proposition}
Combining this proposition with Theorem \ref{jor-bas}, we get:
\begin{theorem}
Let $W(T)$ be a Wick algebra with the braided operator $T$, $-1\le T\le 1$.
Then the kernel of the
Fock representation is generated as a $*$-ideal by $\ker (1+T)$.
\end{theorem}
This theorem implies that, for $q_{ij}$-CCR, $\vert q_{ij}\vert<1$, the Fock
representation is faithful. For the TCCR Wick algebra, and for $q_{ij}$-CCR,
the kernels of the Fock representations are generated by the families
$a_ja_i-\mu a_ia_j$, $i<j$, and $a_ja_i-q_{ij}a_ia_j$, $i<j$, respectively;
and hence the Fock representations of quotients of these algebras by
the $*$-ideals generated by these families are faithful.

\section{Universal bounded representation}\label{jor-stab}
In this section we discuss universal enveloping $C^*$-algebras for
$q_{ij}$-CCR and Wick TCCR.

Let us recall that the universal $C^*$-algebra
for a certain $*$-algebra $\mathcal{A}$ is also called
the universal bounded representation.
It is
the $C^*$-algebra $\mathbf{A}$ with natural homomorphism
$\psi\colon\mathcal{A}\rightarrow\mathbf{A}$ such that, for any
homomorphism $\varphi\colon\mathcal{A}\rightarrow B$, where $B$ is a
$C^*$-algebra, there exists
a unique homomorphism
$\theta\colon\mathbf{A}\rightarrow B$
satisfying $\theta\circ\psi=\varphi$. It can be
obtained by the completion of $\mathcal{A}/J$ with the following
$C^*$-seminorm on $\mathcal{A}$:
\[
\Vert a\Vert=\sup_{\pi}\Vert \pi(a)\Vert,
\]
where $\sup$ is taken over all bounded representations of
$\mathcal{A}$, and $J$ is the kernel of this seminorm. Obviously
this process requires that
$\sup_{\pi}\Vert \pi(a)\Vert<\infty$ for any $a\in\mathcal{A}$.
Note that for our examples this condition is satisfied.

The universal bounded representation for $q_{ij}$-CCR was studied
in \cite{jor-dn,jor-jsw2}.
The following proposition follows from the main result of paper \cite{jor-jsw2}.
\begin{proposition}
Let $\mathbf{A}_{\{q_{ij}\}}$ be the universal enveloping
$C^*$-algebra for $q_{ij}$-CCR,
$\vert q_{ij}\vert<\sqrt{2}-1$. Then there exists the natural isomorphism
\[
\mathbf{A}_{\{q_{ij}\}}\cong\mathbf{A}_0,
\]
where $\mathbf{A}_{0}$ is a $C^*$-algebra generated by the isometries
$ s_i$, $i=1,\ldots,d$, satisfying
\[
s_i^* s_j=0,\ i\ne j
\]
i.e., isomorphism with the Cuntz-Toeplitz algebra.
\end{proposition}
This implies that the
Fock representation of $\mathbf{A}_{\{q_{ij}\}}$ is faithful.

Let us consider the $\mathbf{A}_{\{q_{ij}\}}$, $\vert q_{ij}\vert=1$, for
any $i\ne j$ and $q_{ii}:=q_i$, $\vert q_i\vert<1$
(i.e., unimodular off-diagonal terms). In this case, we do not have
stability on the whole set of parameters (see \cite{jor-dpl}).
\begin{proposition}
If for any $i\ne j$ we have $\vert q_{ij}\vert =1$, then
$\mathbf{A}_{\{q_{ij}\}}$
is isomorphic to the $C^*$-algebra $\mathbf{A}_{0,\{q_{i}\}}$
generated by isometries $\{s_i,\ i=1,\ldots,d\}$ satisfying
\[
s_i^*s_j=q_{ij}s_js_i^*,\ s_js_i=q_{ij}s_is_j,\ i\ne j,
\]
and the Fock representation of $\mathbf{A}_{\{q_{ij}\}}$ is faithful.
\end{proposition}
Finally for the universal $C^*$-algebra $\mathbf{A}_{\mu}$ for the
Wick TCCR, we have the isomorphism $\mathbf{A}_{\mu}\cong
\mathbf{A}_0$ for any $-1<\mu<1$, where the $C^*$-algebra
$\mathbf{A}_0$ is generated by the partial isometries
$\{s_i,\ i=1,\ldots,d\}$ satisfying the relations
\[
s_i^*s_j=\delta_{ij}\left(1-\sum_{k<i}s_k s_k^*\right),\quad i,j=1,\ldots,d.
\]
The Fock representation of $\mathbf{A}_{\mu}$ is faithful also (see
\cite{jor-prsm}).\medskip

\noindent
\textsc{Acknowledgements.}
P. J. was partially supported by the NSF under grants
DMS-9700130 and INT-9722779.

\end{article}

\providecommand{\bysame}{\leavevmode\hbox to3em{\hrulefill}\thinspace}


\begin{thebibliography}{99}

\bibitem{jor-bied} L.C. Biedenharn, 
\textit{The quantum group $\mathrm{SU}_q(2)$ and a
$q$-analogue of the boson operators}, {J. Phys. A}
\textbf{22} (1989), L873--L878.

\bibitem{jor-burkl} I.M. Burban and A.U. Klimyk, \textit{On spectral properties of
$q$-oscillator operators}, {Lett. Math. Phys.} \textbf{29} (1993),
13--18.

\bibitem{jor-fiv} D.I. Fivel, \textit{Interpolation between Fermi and Bose
statistics using generalized commutators}, {Phys. Rev. Lett.}
\textbf{65} (1990), 3361--3364.

\bibitem{jor-jps} P.E.T. J{\o}rgensen, D.P. Proskurin and
Yu.\ S. Samo\u\i{}lenko, \textit{The kernel of Fock representations of Wick
algebras with braided operator of coefficients}, accepted for
publication in {Pacific J. Math.}, math-ph/0001011.

\bibitem{jor-mcf} A.J. Macfarlane, \textit{On $q$-analogues of the quantum
harmonic oscillator and the quantum group $\mathrm{SU}(2)_q$}, {J. Phys.
A} \textbf{22} (1989), 4581--4588.

\bibitem{jor-pw} W. Pusz and S.L. Woronowicz, \textit{Twisted second quantization},
{Rep. Math. Phys.} \textbf{27} (1989), 251--263.

\bibitem{jor-rfw1} R.F. Werner, \textit{The free quon gas suffers Gibbs' paradox},
{Phys. Rev. D (3)} \textbf{48} (1993), 2929--2934.

\bibitem{jor-bsp} M. Bo\.zejko and R. Speicher, \textit{An example of a
generalized Brownian motion}, {Commun. Math. Phys.}
\textbf{137} (1991), 519--531.

\bibitem{jor-dn} K. Dykema and A. Nica, \textit{On the Fock representation of
the $q$-commutation relations}, {J. Reine Angew. Math.}
\textbf{440} (1993), 201--212.

\bibitem{jor-jsw2} P.E.T. J{\o}rgensen, L.M. Schmitt and R.F. Werner,
\textit{$q$-canonical commutation relations and stability of the Cuntz
algebra}, {Pacific J. Math.} \textbf{165} (1994), 131--151.

\bibitem{jor-jsw}
\bysame, \textit{Positive representations of general commutation relations
allowing Wick ordering}, {J. Funct. Anal.} \textbf{134}
(1995), 33--99.

\bibitem{jor-bs}
M. Bo\.zejko and R. Speicher, \textit{Completely positive maps on
Coxeter groups, deformed commutation relations, and operator
spaces}, {Math. Ann.} \textbf{300} (1994), 97--120.

\bibitem{jor-mral} W. Marcinek and R. Ralowski, \textit{On Wick algebras with
braid relations}, {J. Math. Phys.} \textbf{36} (1995), 2803--2820.

\bibitem{jor-green} O.W. Greenberg, \textit{Particles with small violations of
Fermi or Bose statistics}, {Phys. Rev. D (3)} \textbf{43} (1991),
4111--4129.

\bibitem{jor-marc} W. Marcinek, \textit{On commutation relations for quons},
{Rep. Math. Phys.} {\bf 41} (1998), 155--172.

\bibitem{jor-jswc} P.E.T. J{\o}rgensen and R.F. Werner, \textit{Coherent states of
the $q$-canonical commutation relations}, {Commun. Math. Phys.}
\textbf{164} (1994), 455--471.

\bibitem{jor-cun} J. Cuntz, \textit{Simple $C^*$-algebras generated by isometries},
{Commun. Math. Phys.} \textbf{57} (1977), 173--185.

\bibitem{jor-dpl} 
D. Proskurin, \textit{Stability of a special class of
$q_{ij}$-CCR and extensions of higher-dimensional noncommutative tori},
to appear in {Lett. Math. Phys.}

\bibitem{jor-prsm} D. Proskurin and Yu. Samoilenko, \textit{Stability of a $C^*$
algebra associated with the TCCR}, submitted to {Algebras and
Representation Theory}.

\bibitem{jor-osam1} {V. Ostrovsky\u\i{} and Yu. Samo\u\i{}lenko},
\textit{Introduction
to the Theory of Representations of Finitely Presented $*$-Algebras, I:
Representations by bounded operators}, {The Gordon and Breach
Publishing Group, London}, 1999.

\bibitem{jor-proid}
D. P. Proskurin, \textit{Homogeneous ideals in Wick $*$-algebras},
{Proc. Amer. Math. Soc.} \textbf{126} (1998),
3371--3376.
\end{thebibliography}
\end{document}